\documentclass[12pt]{amsart}
\usepackage{amsfonts}
\usepackage[dvips]{graphicx}
\usepackage{amsmath,amsthm}

\newtheorem{theo+}    {Theorem}      [section]
\newtheorem{prop+}  [theo+]  {Proposition}
\newtheorem{coro+}  [theo+]  {Corollary}
\newtheorem{lemm+}  [theo+]  {Lemma}
\newtheorem{deep+}  [theo+]  {Deep Result}
\newtheorem{fact+}  [theo+]  {Fact}
\theoremstyle{definition}
\newtheorem{exam+}  [theo+]  {Example}
\newtheorem{rema+}  [theo+]  {Remark}
\newtheorem{defi+}  [theo+]  {Definition}
\newtheorem{xca+}[theo+]{Exercise}

\numberwithin{equation}{section}

\def\beqn{\begin{equation}}
\def\eeqn{\end{equation}}

\def\epf{\qed \enddemo}

\def\Del{\Delta}

\def\fg{\mathfrak g}
\def\fk{\mathfrak k}

\def\fp{\mathfrak p}

\def\a{\alpha}

\def\Claminv2{|C(\Lambda)|^{-2}}
\def\ga{\gamma}

\def\Ome{\Omega}

\def\Aa2D{A^{\a,2}(D)}
\def\bAa2D{\overline{A^{\a,2}(D)}}
\def\Ab2D{A^{\beta,2}(D)}
\def\bAb2D{\overline{A^{\beta,2}(D)}}

\def\ainnerp#1#2{\langle#1, #2\rangle}

\def\Norm#1_#2{\Vert#1\Vert_{#2}}

\def\2pd#1#2{\frac{\partial^2 #1}{\partial #2^2}}
\def\p11d#1#2#3{\frac{\partial^2 #1}{  \partial #2\partial #3  }}

\def\a{\alpha}

\def\ga{\gamma}
\def\Claminv2{|C(\Lambda)|^{-2}}

\def\Del{\Delta}

\def\det{\operatorname{det}}

\def\Tr{\operatorname{Tr}}

\def\Aa2D{A^{\a,2}(D)}
\def\bAa2D{\overline{A^{\a,2}(D)}}
\def\Ab2D{A^{\beta,2}(D)}
\def\bAb2D{\overline{A^{\beta,2}(D)}}

\def\m{\underline{\bold m}}
\def\ub1#1{\underline{\bold 1^{#1}}}

\def\det{\operatorname{det}}

\def\zae{\"a{}}

\def\bpf{\begin{proof}}
\def\epf{\end{proof}}
\def\beq{\begin{equation}}
\def\eeq{\end{equation}}

\begin{document}

\title[Nearly Holomorphic Functions and
Relative Discrete Series]{Nearly Holomorphic Functions and
Relative Discrete Series of Weighted $L^2$-Spaces on
Bounded Symmetric Domains}

\author{Genkai Zhang}
\address{Department of Mathematics, University of Karlstad,
S-651 88 Karlstad, Sweden}
\address{(current address) Department of Mathematics,
Chalmers University of Technology and G\"o{}teborg University,
S-412 96 G\"oteborg, Sweden}
\email{genkai@math.chalmers.se}


\begin{abstract}
Let $\Ome=G/K$ be a bounded symmetric domain in a complex vector
space $V$ with the Lebesgue measure $dm(z)$
and the Bergman reproducing kernel
 $h(z, w)^{-p}$. Let $d\mu_{\a}(z)=h(z, \bar z)^{\a}dm(z)$,
$\a>-1$, be the weighted measure on $\Ome$.
The group $G$  acts unitarily on the space $L^2(\Ome, \mu_{\a})$
 via change
of variables together with a multiplier.
We consider the discrete parts, also called
the relative discrete series, in the
irreducible decomposition
of the   $L^2$-space.
Let $\bar D=B(z, \bar z)\partial $ be the  invariant Cauchy-Riemann
operator. 
 We realize 
the relative discrete series as the kernels
of the power $\bar D^{m+1}$ of the invariant Cauchy-Riemann operator
$\bar D$ and thus as nearly holomorphic functions
in the sense of Shimura.
We prove that, roughly speaking,
 the operators $\bar  D^m$ 
are intertwining operators from 
the relative discrete series into 
the standard modules of holomorphic discrete series (as
Bergman spaces of  vector-valued holomorphic functions on $\Ome$).
\end{abstract}

\maketitle

\baselineskip 1.40pc

\begin{section}{Introduction}
\end{section}

Let $\Ome$ be a bounded symmetric domain in a complex vector
space $V$ with the Lebesgue measure $dm(z)$.
The Bergman reproducing kernel is up to a constant
$h(z, \bar w)^{-p}$, where $h(z, \bar w)$ is an irreducible
polynomial holomorphic in $z$ and antiholomorphic in $w$. We consider
the weighted measure $d\mu_{\a}(z)=h(z, \bar z)^{\a}dm(z)$ for 
$\a>-1$ and corresponding $L^2$-space $L^2(\Ome, \mu_{\a})$
on $\Ome$. The group $G$ of biholomorphic mappings
of $\Ome$ acts unitarily on the $L^2$-space via change
of variables together with a multiplier, and the weighted
Bergman space is then an irreducible invariant subspace. 
The irreducible  decomposition of the 
$L^2$-space under the $G$-action
has been given by Shimeno \cite{Shimeno}. It is proved
there abstractly (via
identifying the infinitesimal characters)
that all the discrete parts (called relative
discrete series) appearing in the decomposition
are holomorphic discrete series. In this paper we 
consider their explicit realization.

To illustrate our main results
we consider the case of the unit disk.
The Bergman reproducing kernel is  $(1-z\bar w)^{-2}$, and the
weighted measure in question is $d\mu_{\a}(z)=(1-|z|^2)^\a dm(z)$.
The group $G=SU(1, 1)$ acts unitarily on $L^2(D, \mu_\a)$ via a projective
representation 
$$\pi_{\nu}(g)f(z)= f(g^{-1}z)(cz +d)^{-\nu}, \quad
 g^{-1}=\begin{pmatrix}
 a& b\\c&d\end{pmatrix}
$$
where $\nu=\a+2$.
 To study
the relative discrete series we introduce
the invariant Cauchy-Riemann operator $\bar D=(1-|z|^2)^2\bar \partial $.
The operator $\bar D$ intertwines the action $\pi_{\nu}$
with the action $\pi_{\nu-2}$, which can be proved by direct calculation.
 The kernel $\ker \bar D$ of $\bar D$ on the
weighted $L^2$-space is the weighted Bergman space 
$L^2_a(\Ome, \mu_\a)$ 
of holomorphic functions, which
gives one of the relative discrete series. It is naturally
to expect that the kernel $\ker \bar D^{m+1}$ of the iterate
of $\bar D$ will give us the other relative discrete series.
The functions that are in the kernel $\ker D^{m+1}$ can be written
as polynomial of $q(z)=\frac{\bar z}{1-|z|^2}$ of degree  $\le m$
with coefficients
being holomorphic functions.
Those functions, following
Shimura, are called \textit{nearly holomorphic functions}.
 The function
$q(z)$ actually is the holomorphic differential of the K\zae{}hler potential
$\log (1-|z|^2)^{-2}
$. Indeed $q(z)=\frac 12 \partial_z \log (1-|z|^2)^{-2}
$. Moreover it has a Jordan theoretic meaning as the \textit{quasi-inverse}
of $\bar z$ with respect to $z$ in  $\bar\mathbb C$ with 
the Jordan triple product $\{\bar u  z \bar v\}= 2\bar u z \bar  v$.
The key result is
that each power $q(z)^m=\frac{\bar z^m}{(1-|z|^2)^m}$, for
$0\le m< \frac{\a+1}2$ generates a relative discrete series.
 Denote corresponding
 the relative discrete
  series by
$A^{2, \a}_m$.
Then the operator $\bar D^m$
is an intertwining operator from
$A^{2, \a}_m$ into the weighted Bergman space
in $L^2(\Ome, \mu_{\a-2m})$, namely $L^2_a(\Ome, \mu_{a-2m})$.
Moreover  all relative discrete series are obtained in
this way.

When $\Ome=G/K $ is a general bounded symmetric domain 
the corresponding function $q(z)$, defined as the differential
of the K\"ahler potential, can indeed be expressed in term of
quasi-inverse in the Jordan triple $V$. See Proposition 3.1.
Let $\bar D$ be the invariant Cauchy Riemann operator.
Then it is proved in \cite{pz-cr} that
the iterate $\bar D^m$ maps a function
on $\Ome$ to a function with value in the symmetric
subtensor space $S_m(V)$ of $\otimes^m V$.
Decompose $S_m(V)$ into irreducible subspaces under
$K$.
 Let
$\m$ be the signature of an irreducible
subspace and $\bar \Del_{\m}$ the highest weight
vector in that space, considered as a polynomial function
on $V^\prime$. Now the function $q(z)$ is a $V^\prime$-valued
function on $\Ome$, thus $\bar\Del_{\m}(q(z))$ is a scalar-valued
function on $\Ome$. We prove that $\bar\Del_{\m}(q(z))$ is
in the space $L^2(\Ome, \mu_{\a})$ when $\m$
satisfies certain condition; see Proposition 4.1. We further
prove that it generates an irreducible subspace, namely
a relative discrete series, and is the highest weight vector, and that
 the operators $\bar D^m$ 
are intertwining operators from the relative
discrete series onto the weighted
Bergman space of holomorphic functions
with values in the irreducible subspace
$S_{\m}(V)$ of the symmetric tensor $S_m(V)$, the
later being a  standard
module of holomorphic discrete series.
We thus realize 
the relative discrete series in the kernel
of the power $\bar D^{m+1}$ and  as nearly holomorphic functions
in the sense of Shimura (\cite{Shimura-annmath-86}
and \cite{Shimura-mathann-87}). 

Finally in the last section we consider as an example  the unit ball in 
$\mathbb C^n$. We calculate  directly, via
the adjoint operator $\bar D^\ast$,
the highest weight vector in the relative discrete
series. The realization of the relative discrete series
has  also been studied in \cite{pz-cr}.

Our results explain geometrically why
the relative discrete series are equivalent to
the weighted Bergman spaces with values in
\textit{symmetric tensor space} of the tangent space. Moreover,
since the highest weight vectors are quite explicitly
given we understand better the analytic nature of the functions
in the discrete series. 
 We hope that our result will
be helpful in understanding the $L^p$-spectral properties
of the irreducible decomposition, for example,
the $L^p$-boundedness of the orthogonal projection into
the relative discrete series.

\subsection*{Acknowledgments. } 
I would like to thank  Jaak Peetre and Harald Upmeier 
for  some illuminating discussions. 
I would also like to thank the Erwin Schr\"o{}dinger Institute
for mathematical physics, Vienna, for providing a
 stimulating environment.

\begin{section}{Invariant
Cauchy-Riemann  Operator $\bar D$
and Nearly Holomorphic Functions on K\zae{}hler manifolds
}
\end{section}

We recall in this section briefly some preliminary results
on invariant Cauchy-Riemann operators
and nearly holomorphic functions on K\zae{}hler manifolds;
see \cite{Shimura-annmath-86}, \cite{Englis-Peetre}, 
 \cite{gz-shimura}, and \cite{gz-invdiff}.

Let  $\Ome$
be a  K\zae{}hler manifold with the 
K\zae{}hler metric locally  given by the matrix  $(h_{i\bar
j})$, with $h_{i\bar j}=\frac{\partial^2 \Psi}{\partial z_i\partial
\bar z_j}$ and a potential  $\Psi$. Let $T^{(1, 0)}$
be its holomorphic tangent bundle.
Let $W$ be a Hermitian vector bundle
over $\Ome$, and $C^\infty(\Ome, W)$
its smooth sections.
The invariant Cauchy-Riemann
operator $\bar D$ locally  defined as follows.
If $f=\sum_{\a}f_\a e_\a$ is any section of $W$, then
$$
\bar D f=\sum_{\a, i, j} h^{\bar \jmath i} \frac{\partial f_{\a}} {\partial
  {\bar z^j}}\partial_i \otimes e_{\a}.
$$
It maps $f\in C^\infty(\Ome, W)$ to $\bar Df\in C^\infty(\Ome,
T^{(1, 0)}\otimes W)$.
Denote  $S_m(T^{(1, 0)})
$ the symmetric tensor subbundle
of $\otimes^mT^{(1, 0)}$.
We recall some known properties of the operator $\bar D$. See \cite{pz-cr}.

\begin{lemm+} The following assertions hold.
\begin{description}
\item[(1)]
The operator $\bar D$ is an intertwining operator: If $g$
is a biholomorphic mapping of $\Ome$, then
\begin{equation}\label{cov-pro}
\bar D(g_W f)=((dg)^{-1}\otimes g_W) Df,
\end{equation}
where $g_W$ is the induced action  of $g$ on sections of $W$ and $dg(z)$:
$T^{(1, 0)}_z\mapsto T^{(1, 0)}_{gz}$ is the differential of $g$.
\item[(2)]The iterate $\bar D^m$ of $\bar D$
 maps $C^{\infty}(\Ome, W)$
to $C^{\infty}(\Ome, W\otimes S_m(T^{(1, 0)}))$. 
\end{description}
\end{lemm+}

For our later purpose we can assume that $\Ome$ is some
domain in a vector space $V$ with coordinates $\{z_j\}$,
and that all the bundles are trivial. The space $T^{(1, 0)}_z$
will be identified with $V$.
So let $W$ be a vector space and we will consider the 
$C^\infty(\Ome, W)$
 of $W$-valued $C^\infty$-functions on
$\Ome$.

Let 
$$
q(z)=\partial \Psi =\sum_j \frac{\partial\Psi}{\partial z_j} dz_j.
$$
Here
$\Psi$ is the K\"a{}hler potential, and $\{dz_j\}$ is the dual basis
 for the holomorphic
cotangent space $V^\prime$. Thus $q(z)$ is a function
with values in $V^\prime$.
 Following Shimura \cite{Shimura-annmath-86}
we call a  $W$-valued function $f\in C^\infty(\Ome, W)$
nearly holomorphic if $f$ is a polynomial
of $q(z)$ with holomorphic coefficients.
We denote $\mathcal N_m$  the space of scalar-valued
nearly holomorphic
functions that are polynomial of degree $\le m$, namely
those functions $f(z)=\sum_{|\underline{\beta}|\le m}c_{\underline{\beta}}(z) q(z)^{\underline{\beta}
}
$
where $c_{\underline{\beta}}(z)$ are holomorphic 
functions.
 
We denote 
 $\text{Id}$ 
the identity tensor in the tensor
product $V\otimes V^{\prime}$.  By the direct calculation
 we have
\begin{equation}\label{dn}
\bar D q(z)= \text{Id};
\end{equation}
see \cite{gz-shimura}. We generalize this  formula  as follows; the proof
of it is quite straightforward and we omit it.
\begin{lemm+} We have the following differentiation formula
\begin{equation}\label{dnm}
\bar D^m (\otimes^m q(z))= m!\text{Id},
\end{equation}
where $\text{Id}$  in the right hand denotes
the identity tensor in the tensor product
$(S_mV)\otimes S_m (V^{\prime})=
(S_m V)\otimes(S_m V)^{\prime}$.
\end{lemm+}

\begin{rema+} The formula (\ref{dn}) was observed earlier by Shimura
\cite{Shimura-annmath-86} and Peetre \cite{Peetre-cr-manu}; in the
later paper explicit formulas were given 
for the Laplace operators
on weighted $L^2$-spaces on bounded symmetric domains, where the
function
$q(z)$ also appears.
\end{rema+} 

\begin{exam+} We consider the case of the unit disk. 
The operator $\bar D=(1-|z|^2)^2\bar \partial$
The function $q(z)$ is $\frac{\bar z}{1-|z|^2}$
(or  exactly it is $\frac{\bar z}{1-|z|^2}dz$).
The 
above formula amounts to 
$$
\bar D^m
(\frac{\bar z}{1-|z|^2})^m
=m!,
$$
which can be proved by direct calculations.
It can also be proved by using the formula
$$
\bar D^m= (1-|z|^2)^{m+1}
(\frac{\partial}{\partial \bar z})^m
(1-|z|^2)^{m-1},
$$
see \cite{gz-shimura}. Indeed, 
\begin{equation*}
\begin{split}
&\quad \bar D^m(\frac
{\bar z}{1-|z|^2})^m\\
&= (1-|z|^2)^{m+1}
(\frac{\partial}{\bar \partial z})^m
(1-|z|^2)^{m-1}(\frac{\bar z}{1-|z|^2})^m\\
&= (1-|z|^2)^{m+1}
(\frac{\partial}{\bar \partial z})^m
\frac{\bar z^m}{1-|z|^2}\\
&= (1-|z|^2)^{m+1}\sum_{l=0}^m\binom{m}{l} m(m-1)\cdots(m-l+1)\bar z^{m-l}
\frac{(m-l)! z^{m-l}}{(1-|z|^2)^{m-l+1}}\\
&=m!\sum_{l=0}^m\binom{m}{l}(z \bar z)^{m-l}(1-|z|^2)^{l}
\\
&=m!
\end{split}
\end{equation*}
The calculations are somewhat combinatorially intriguing.
\end{exam+}

Using the
above result we get immediately
the following characterization of nearly holomorphic functions.
This is proved in \cite{Shimura-annmath-86}, Proposition 2.4,
for classical domains. It can be proved for all
K\zae{}hler manifolds via the same methods.

\begin{lemm+}\label{ker-d-m}Consider the operator $\bar D^{m+1}$
on the space $C^{\infty}(D)$
of  $C^\infty$-functions on $D$. Then
$$
Ker \bar D^{m+1} =\mathcal N_m.
$$  
\end{lemm+}

We recall the identification of polynomial functions
with symmetric tensors. This will clarify conceptually some
calculations in the next section.
There is a  pairing 
$$
(\phi, \psi)\in S_m(V)\times S_m(V^\prime)\mapsto [\phi, \psi]\in \mathbb C,
$$
between the symmetric tensor spaces  $S_m(V)$
and $S_m(V^\prime)$,
via the natural
pairing between 
$\otimes^m V$ and 
$\otimes^m V^\prime$.
Now for each 
element $\phi$ in the symmetric tensor space $S_m(V)$
there corresponds a homogeneous polynomial function
of degree $m$ on the
space $V^\prime$,
also denoted by $\phi$,
such that
\begin{equation}\label{pol-sym-iden}
[\phi, v^\prime\otimes v^\prime\otimes \cdots v^\prime]=\phi(v^\prime)
\end{equation}
for any $v^\prime\in V^\prime$.

Using this convention we see that a function $f\in C^\infty(\Ome)$
is in $\mathcal N_m$
if and only if there exist holomorphic functions $g_k$
with values in the tensor product $ S_k(V)$, $k=0, 1, \dots, m,$
such that
\begin{equation}\label{exp-nhf}
f(z)=\sum_{k=0}^m g_k(q(z)).
\end{equation}

\begin{section}{Nearly Holomorphic Functions on Bounded Symmetric Domains}
\end{section}

In this section we assume that $\Ome=G/K$ is a bounded symmetric
domain of rank $r$ in a complex vector space $V$. Here $G$ is the identity component of
the group of biholomorphic mappings of $\Ome$ and $K$
is the isotropy group at $0\in V$. Let $\mathfrak g =\mathfrak k +\mathfrak p$
be the Cartan decomposition of $\mathfrak g$.
The space $V$ has  a Jordan triple structure
so that the space $\mathfrak p$ is explicitly described; see 
\cite{Loos-bsd}, whose notation and results will be incorporated
here. So let $Q(z): \bar V\to V$ be the quadratic operator. The
$$
\mathfrak p=\{\xi_v=v-Q(z)\bar v\}
$$
viewed as holomorphic vector fields on $\Ome$. Let
$D(z, \bar v)w=\{z\bar v w\}=(Q(z+w)-Q(z)-Q(w))\bar v$
be the Jordan triple product.
We normalize the $K$-invariant
Hermitian inner product $\langle z,w \rangle$
on $V$ so that a minimal tripotent has norm
$1$. This can also be calculated by
\begin{equation}\label{normal-V}
\langle z,w \rangle=\frac 1p \Tr D(z, \bar w)
\end{equation}
where $p$ is an integer called the genus of $\Ome$. We identify then
the vector space $V^\prime $ with $\bar V$ via this scalar product.

Let $dm(z)$ be the corresponding Lebesgue
measure on $V$. The Bergman reproducing kernel on $D$
is the $ch(z, w)^{-p}$ for some positive constant $c$.
Let
$$
B(z, \bar w)=I-D(z, \bar w)+Q(z)Q(\bar w)
$$
the Bergman operator. $B(z, \bar w)$ is holomorphic
in the first argument and anti-holomorphic in the second.
(We write $B(z, \bar w)$ instead of $B(z, w)$ as in \cite{Loos-bsd}
in order to differ it from $B(\bar z, w)$ which is acting on the space 
$\bar V$.)
 The Bergman metric
at $z\in \Ome$
defined by  the metric   $\partial_j\bar \partial_k
\log h(z, \bar z)^{-p}$ on $\Ome$
is then
$$
p\langle B(z, \bar z)^{-1}z, w\rangle;
$$
and
$$
\det B(z, \bar z)= h(z, \bar z)^p
$$
See \cite{Loos-bsd}. 
 For some computational
convenience we will choose and fix
the metric on $\Ome$ to be
\begin{equation}
\boxed{\langle B(z, \bar z)^{-1}z, w\rangle}.
\end{equation}

The invariant Cauchy-Riemann operator is
$$
\bar D =B(z,\bar z)\bar\partial,
$$
and the $N$-function defined in the previous section is
now (with a normalizing constant)
\begin{equation}\label{N-bsd}
\boxed{q(z)=\frac 1p \partial \log \det B(z, \bar z)^{-1}}
\end{equation}

We shall find an explicit formula for the function $q(z)$
on  $\Ome$. Recall first the notion of 
\textit{quasi-inverse}  in the Jordan triple $V$;
see \cite{Loos-bsd}. Let $z\in V$
and $\bar w\in \bar V$.
The element $z$ is called quasi-invertible with respect to
$w$ if $B(z, \bar w)$ is invertible and its quasi-inverse
is given by
$$
z^{\bar w}=B(z, \bar w)^{-1}(z-Q(z)\bar w).
$$
Similarly we define the quasi-inverse of an element $\bar z\in \bar V$
with respect to $w\in V$.
\begin{prop+} The function $q(z)$ on $\Ome$ is
given by
$$
q(z)=\bar z^{z}
=B(\bar z,  z)^{-1}(\bar z-Q(\bar z) z)
$$
\end{prop+}

\begin{proof} For some
computational
convenience  we  consider, instead of the holomorphic differential in
 (\ref{N-bsd}),
 the anti-holomorphic differential
$$
\frac 1p\bar \partial\log  \det B(z, \bar z)^{-1}=
-\frac 1p\bar \partial\log  \det B(z,\bar  z).
$$ 
 Let $\bar v\in \bar  V$. By the 
definition of $B$-operator
we have
\begin{equation}
\begin{split}
&\quad\, B(z, \bar z+t\bar v)\\
&=1-D(z, \bar z+t\bar v)+Q(z)Q(\bar z+t\bar v)\\
&=1-D(z, \bar z)+Q(z)Q(\bar z)+ t (-D(z, \bar v)+Q(z)Q(\bar z, \bar v)
+t^2 Q(z)Q(\bar v)\\
&=B(z, \bar z)\left(I + tB(z, \bar z)^{-1}(-D(z, \bar v)
+Q(\bar z, \bar v)) +t^2B(z, \bar z)^{-1}Q(z)Q(\bar v)\right).
\end{split}
\end{equation}
Thus
the first order term in $t$ in $\log \det B(z, \bar z +t \bar v)$ is 
\begin{equation}\label{1-order}
 \Tr (B(z, z)^{-1}(-D(z, \bar v)+ Q(z)Q(\bar z, \bar v)).
\end{equation}

We recall a formula  in \cite{Loos-bsd} (see (JP30))
$$
B(z, \bar z)D(z^{\bar z}, v)=D(z, \bar v)- Q(z)Q(\bar z, \bar v).
$$
Therefore (\ref{1-order}) is
$$
-\Tr D(z^{\bar z}, v)=-p\langle z^{\bar z}, v\rangle
$$
by the formula (\ref{normal-V}).
Summarizing we find
$$
\frac 1p
\bar\partial_{v}\log  \det B(z, \bar z)^{-1}=
\langle z^{\bar z}, v\rangle,
$$
which is the desired formula.
\end{proof}

Now the group $K$ acts on $\Ome$ and keeps the function $h(z, \bar z)$-invariant.
Thus we get, in view of the formula (\ref{N-bsd}),
\begin{equation}\label{k-N}
q(kz)=(k^{-1})^\prime q(z),
\end{equation}
where $(k^{-1})^\prime $ on $q(z)\in V^\prime$
is the dual of $k^{-1}$ on $V$.

 In particular, since 
the
function $q(z)$ is a $V^\prime$-valued function on $\Ome$, we have, for any
homogeneous polynomial function $f$  on $V^\prime$,
a scalar-valued function $f(q(z))$.
The following lemma  then 
follows from (\ref{k-N}) and the $K$-invariance of the pairing
between $S_m(V^\prime)$ and $S_m(V)$.

\begin{lemm+}\label{k-intert} 
The map
 $$ v\in S_m(V)\mapsto v(q(z))= [v, \otimes^m q(z)]
$$
is an invertible $K$-intertwining operator between the $K$-action
on $S_m(V)$ and its regular action on  functions on $\Ome$.
\end{lemm+}

We recall now the decomposition of $S_m(V)$ under $K$.
To state the result we fix
some notation.
The complexification $\fg^{\mathbb C}$
of the  Lie algebra  $\fg$ has
a decomposition $\fg^{\mathbb C}=\fp^{+}+
\fk^{\mathbb C}+\fp^{-}$, with $\fk^{\mathbb C}$
the complexification of
 the Lie algebra $\fk$ of $K$ and $\fp^{+}=V$.
Let $\{e_1, \dots, e_r\}$ be a frame of tripotents in  $V$.
Fix  an Cartan subalgebra of
of $\mathfrak k^{\mathbb C}$, and let $\ga_1> \cdots >\ga_r$
be the Harish-Chandra strongly
roots so that $e_1, \dots, e_r$
are the corresponding root
vectors. The ordering of the roots of $\fg^{\mathbb C}$
is so that $\fp^{+}$ is the sum of positive non-compact root vectors. We shall
then  speak of \textit{highest weight modules} of $\fg^{\mathbb C}$
with respect to this ordering.

\begin{lemm+}\label{Hua}
(\cite{Hua}, \cite{Schmid} and \cite{FK}) The space 
$S_m(V)$ (respectively $S_m(V^\prime)$) under $K$
 is decomposed into irreducible subspaces 
with multiplicity one as
$$
S_m(V)=\sum_{\m}S_{\m}(V), \qquad (\text{resp.}\,  S_m(V^\prime)=
\sum_{\m}S_{\m}(V^\prime))
$$
where each $S_{\m}(V)$ (resp. $S_{\m}(V^\prime)$) 
is of highest weight $\m=m_1\ga_1 +\cdots +m_r
\ga_r$ (resp. lowest weight $-(m_1\ga_1+\cdots + m_r\ga_r)$)
with $ m_1\ge m_2\ge\cdots\ge m_r\ge 0$,
and the summation is over all
 $\m$
with $|\m|=m_1+m_2 +\cdots +m_r =m$.
\end{lemm+}

The highest weight vectors of $S_{\m}(V)$
(respectively  lowest weight vectors of $S_{\m}(V^\prime)$)
 have constructed explicitly;
see \cite{FK} and reference therein. Let $\Del_j$
be the lowest weight vector of the fundamental representation
$\m=\ub1j=\ga_1+\dots+\ga_j$, $j=1, \dots, r$. The polynomial
$\Del=\Del_r$ is the determinant function of the Jordan triple $V$.
Then the lowest weight vector of $S_{\m}(V^\prime)$ is  
\begin{equation}
\Del_{\underline{\bold{m}}}(v)=\Del_1(v)^{m_1-m_2}
 \cdots  \Del_{r-1}(v)^{m_{r-1}-m_r}
\Del_{r}(v)^{m_r},
\end{equation}
viewed as polynomial of $v\in V$. 
Via the natural
pairing between $S_{\m}(V^\prime)$ 
and $S_{\m}(V)$
 we  find that  the highest weight vector
of $S_{\m}(V)$ is  $\bar \Del_{\m}$ and
\begin{equation}\label{con-hwv}
\bar \Del_{\underline{\bold{m}}}(w)=\bar \Del_1(w)^{m_1-m_2}
 \cdots  \bar \Del_{r-1}(w)^{m_{r-1}-m_r}
\bar\Del_{r}(w)^{m_r},
\end{equation}
viewed as polynomial of $w\in V^\prime=\bar V$.

\begin{section}{The Relative Discrete Series of
$L^2(\Ome, \mu_{\a})$}
\end{section}

In this section we find a family of relative discrete series
by constructing some vectors that are in   $L^2$-space and
are highest weight vectors, namely annihilated by the positive vectors
in $\fg^{\mathbb C}$ via the induced action of (\ref{pi-nu}) (see
below).

Let $\alpha >-1$. 
and consider the weighted measure
$$
d\mu_{\a}=h(z,z)^{\a}dm(z).
$$
The group $G$  acts unitarily on the space $L^2(\Ome, d\mu_{\a})$
via
\begin{equation}\label{pi-nu}
\pi_{\nu}(g)f(z)=f(g^{-1}z){J_{g^{-1}}(z)}^{\frac{\nu}{p}}, \quad g\in G.
\end{equation}
where $\nu =\a +p$ and $J_g$ is the Jacobian determinant of $g$.
We denote $L^2_a(\Ome, \mu_\a)$ 
the  weighted Bergman 
space of holomorphic functions in  $L^2(\Ome, \mu_\a)$. 

We introduce now the  weighted Bergman spaces
of vector-valued holomorphic functions that will be used
to realize the relative discrete series in $L^2(\Ome, d\mu_{\a})$.
Fix  a signature $\m$ with $m=m_1+\dots+m_r$.
We denote 
$L^2_a(\Ome, S_{\m}(V), \mu_\a)$
the weighted Bergman space of $S_{\m}(V)$-valued
holomorphic functions such that the following norm is finite
$$
\Vert f\Vert^2
=\int_{\Ome} \ainnerp{  (\otimes^m B(z, \bar z) ^{-1})  f(z)}{ f(z)}\,
d\mu_{\a}(z).
$$
The group $G$ acts unitarily on $L^2_a(\Ome, S_{\m}(V), \mu_\a)$
via
\begin{equation}\label{g-on-berg}
g\in G: f(z)\mapsto      (J_{g^{-1}}(z))^{\frac \nu p}\otimes ^m (dg^{-1}(z))^{-1}  f(g^{-1}z).
\end{equation}
This space is non trivial and  forms
an irreducible representation of $G$
when $\m$ satisfies the following condition:
\begin{equation}\label{con-rds}
\frac{\a+1}2>m_1\ge m_2 \ge \cdots \ge m_r \ge 0.
\end{equation}
This follows directly from Theorem 6.6 in \cite{Kn-book}; see
also \cite{Shimeno}. (We note here that non-triviality of
the space can also be proved directly
by expressing the inverse $B(z, \bar z)^{-1}$ of the Bergman
operator via the quasi-inverse developed in \cite{Loos-bsd}, quite
similar 
to the proof of Proposition 4.1 below. However we will not go into the 
details here.)

Our first result is a construction of certain vectors
in $L^2(\Ome, \mu_{\a})$. 

\begin{prop+}\label{inL2} Suppose $\m$ satisfies the condition 
\ref{con-rds}. Then 
the functions $\bar{\Del_{\m}}(q(z))$
is in $L^2(\Ome, \mu_{\a})$ and in $\text{Ker} \bar D^{m+1} $.
\end{prop+}

We begin with fundamental representations
$S_{\m}(V)$ with signatures $\m=\ub1j=\ga_1+\dots +\ga_j$
and highest weight vectors $\bar\Del_j$,  $j=1, \dots, r$.

\begin{lemm+}  Then the function $\Del_{j}(q(z))$
is of the form
\begin{equation}\label{del-n-1}
\bar\Del_{j}(q(z))=\frac{P(z, \bar z)}{h(z, \bar z)}
\end{equation}
where $P(z, \bar z)$ is a polynomial in $(z, \bar z)$
of total degree not exceeding $2r$.
In
particular if $j=r$,
\begin{equation}\label{del-n-2}
\bar\Del(q(z))= \frac{\bar \Del(z)}{h(z, \bar z)}.
\end{equation}
\end{lemm+}

\begin{proof} It follows from the Faraut-Koranyi expansion
that
\begin{equation}\label{fk-h}
h(v, \bar w)=\sum_{s=0}^{r}(-1)^rc_sK_{\ub1s}(v,\bar w),
\end{equation}
where $K_{\m}$ is the reproducing kernel of the subspace $\mathcal P^{\m}(V)$
of $\mathcal P(V)$
with signature $\m$ with the Fock-norm $\langle\cdot, \cdot\rangle_{\mathcal F}$ and $c_s$ are positive constants; see  \cite{FK}.
Performing the inner product 
in the Fock space 
of the element
$h(v, \bar w)$ with the function $\Del_{1^j}(v)$ and
using (\ref{fk-h})
we find that
$$
\langle h(\cdot, \bar w), \Del_{j}\rangle_{\mathcal F} =
(-1)^s c_s\Vert\Del_j\Vert_{\mathcal F}^2 \bar \Del_j(\bar w);
$$
namely,
\begin{equation}\label{del-n-3}
\bar\Del_j(\bar w)=
\frac 1{(-1)^s c_{s}}\Vert\Del_j\Vert^2_{\mathcal F} 
\langle h(\cdot, \bar w), \Del_{j}\rangle_{\mathcal F}.
\end{equation}
We take now $\bar w =\bar z^z$.
Recall  \cite{Loos-bsd}, Lemma  7.5, that 
\begin{equation}
\label{loos-h-h}
h(v, \bar z^{z})
=\frac{h(v+z, \bar z)}{h(z, \bar z)}.
\end{equation}
Substituting this into the previous formula we get
\begin{equation}\label{del-n-4}
\bar \Del_j(\bar z^z)=
\frac 1{(-1)^s c_{s} h(z, \bar z)}
\Vert\Del_j\Vert^2_{\mathcal F} 
\langle h(\cdot +z, \bar z), \Del_{j}\rangle_{\mathcal F}
\end{equation}
Since $h(v+z, z)$ is a polynomial in $z$ and $\bar z$ of  degree $2r$,
we see that $\Del_j(\bar z^{z})$
is of the declared form. 

If $j=r$, we can then calculate
$\langle h(\cdot, \bar z), \Del_{r}\rangle_{\mathcal F}
$ further.
Expand $h(v+z,  \bar z)$ 
 again using (\ref{fk-h}). We have
\begin{equation}
\begin{split}
\langle h(\cdot+z, \bar z), \Del_{r}\rangle_{\mathcal F}
&=\sum_{s=0}^{r}(-1)^s c_{s}\langle K_{\ub1s}(\cdot+z, \bar z), \Del_{r}\rangle_{\mathcal F}\\
\\
&=(-1)^r c_{s}
\langle K_{\ub1r}(\cdot+z, \bar z), \Del_{r}\rangle_{\mathcal F},
\end{split}
\end{equation}
because $\Del_r$ is of degree $r$ and it is orthogonal
to those terms of lower  degree.
But 
$$
K_{\ub1r}(z+v, z)=
K_{\ub1r}(v, z) + \dots
$$
where  the rest term is of lower order. Therefore by the same reason
and by the reproducing property,
$$
\langle h(z+\cdot, z), \Del_{r}\rangle_{\mathcal F}
=(-1)^r c_{r}
\langle K_{\ub1r}(\cdot, \bar z), \Del_{r}\rangle_{\mathcal F}
=(-1)^r c_{r}\Vert\Del_r\Vert^2_{\mathcal F}\overline{\Del_r(z)}.
$$
Substituting this into (\ref{del-n-4}) we 
then get (\ref{del-n-2}). 
\end{proof}

\begin{rema+} The norm $\Vert\Del_{\m}\Vert_{\mathcal F}$ is calculated
in \cite{FK}, though we will not need it in the present paper.
\end{rema+}

Recall formula (\ref{con-hwv}) for the highest weight vector
$\bar\Del_{\m}$.
As a corollary we find immediately that 
\begin{coro+}\label{del-n}
 Then the function $\bar \Del_{\m}(q(z))$
is of the form
$$
\bar \Del_{\m}(q(z))=
 \frac{P(z, \bar z)}{h(z, \bar z)^{m_1}}
$$
where $P(z, \bar z)$ is a polynomial in $(z, \bar z)$. 
\end{coro+}

We prove now the Proposition \ref{inL2}.
\begin{proof} We estimate the norm
of $\bar \Del_{\m}(q(z))$ in $L^2(\Ome, \mu_\a)$
by using the above Corollary. 
The polynomial $P(z, z)$ on $\Ome$ is bounded, say $|P(z, z)|\le C$.
We have
$$
\int_{\Ome}|\frac {P(z, z)}{h(z, \bar z)^{m_1}}|^2d\mu_{\a}
\le C \int_{\Ome} h(z, \bar z)^{\a-2m_1}dm(z).
$$
By the condition (\ref{con-rds}) we see that $\a-2m_r > -1$, thus
the above integral is finite (see \cite{FK}), namely the function is in the 
$L^2$-space. That $\bar \Del_{\m}(N(z))$ is in $\text{Ker}D^{m+1}$
follows  directly from Lemma \ref{ker-d-m}. 
\end{proof}

The action $\pi_{\nu}$
 of $G$ on $L^2(\Ome, \mu_\a)$ induces
an action of $\fg^{\mathbb C}$ on the space of $C^\infty$-functions.
We  prove next that the function $\bar\Del_{\m}(q(z))$ is annihilated
by the positive root vectors in $\fg^{\mathbb C }$.
 The element  in
$\fp$,  when viewed as holomorphic vector fields, are of the form
$\xi_v=v-Q(z)\bar v$; thus when acting on $C^\infty$-functions
 on 
$\Ome$ induced from the
regular action of $G$, they are 
$$
(\partial_v 
-\partial_{Q(z)\bar v})f
+(\partial_{\bar v} -\partial_{Q(\bar z) v})f
$$
>From this it follow that the element $v\in \fp^+=V$ acts on 
 $C^\infty$-functions induced from $\pi_{\nu}$ of $G$ is
\begin{equation}\label{pi-nu-v}
\pi_{\nu}(v)f=\partial_v f -\partial_{Q(\bar z) v}f,
\end{equation}
since the infinitesimal
 action of $v\in \fp^+$ is a translation and it will not
contribute in the determinant factor
in  (\ref{pi-nu}).

To study the action of $\fp^+$ on $\bar \Del_{\m}(q(z))$, we calculate first
the differentiation of $q(z)$.

\begin{lemm+} The following differentiation
formulas hold
\begin{equation}\label{d-v-q}
\partial_v q(z) =Q(q(z))v,
 \quad 
\partial_{\bar w}q(z)=B(\bar z, z)^{-1}\bar w.
\end{equation}
In particular if $\bar w=Q(\bar z)v$,
\begin{equation}\label{d-v-q-1}
\partial_{Q(\bar z)v}q(z)=B(\bar z, z)^{-1}Q(\bar z)v=Q(q(z))v,
\end{equation}
and
\begin{equation}\label{van-q}
(\partial_v -\partial_{Q(\bar z)v })q(z)= 0
\end{equation}
\end{lemm+}
\begin{proof} We use  the addition formulas in \cite{Loos-bsd}, Appendix,
for the quasi-inverses. As   special cases we have
\begin{equation}\label{q-i-1}
\bar z^{z+t v} = (\bar z^{z})^{tv}=
B(\bar z^z, t v)^{-1}(\bar z^z-tQ(\bar z^z)v),
\end{equation}
and
\begin{equation}
\label{q-i-2}
(\bar z+ t\bar w)^{z} = \bar z^z + B(\bar z, z)^{-1}B(t\bar w, z^{\bar z})^{-1}
 (t \bar w -Q(t\bar w) z^{\bar z}).
\end{equation}
The first order term in $t$ in (\ref{q-i-1}) is easily seen to be
$$
D(\bar z^z, v)\bar z^z -Q(\bar z^z)v
= Q(\bar z^z)v,
$$
which proves the first formula in
(\ref{d-v-q}). Similarly we  can calculate the first order
term in (\ref{q-i-2}) and prove  the second formula; using this formula
and 
$$
B(\bar z, z)^{-1}Q(\bar z)=Q(\bar z^z) = Q(q(z)),
$$
we get then (\ref{d-v-q-1}). 
\end{proof}
 
We can thus calculate $\pi_{\nu}(v)$ on $\bar\Del_{\m}(q(z))$
by using (\ref{pi-nu-v}). In view of (\ref{van-q}) we have
$$
\pi_\nu(v)\bar\Del_{\m}(q(z))=0.
$$
This, together  with Lemma \ref{k-intert}, implies that
\begin{prop+} The vector $\Del_{\m}(q(z))
$ under the action of $\pi_{\nu}$ of $\fg^{\mathbb C}$ is annihilated
by the positive root vectors.
\end{prop+}

We let $A^{2, \a}_{\m}(\Ome)$
 be the 
subspace of $L^2(\Ome, \mu_{\a})$
generated by the function $\bar\Del_{\m}(q(z))$, for $\m$
given by (4.1).  Thus it is a highest weight representation of $G$.
Now it follows from Lemma 2.2 that
$$
\bar D^m (\bar \Del_{\m}(q(z)))=m! \bar\Del_{\m}
$$
The vector $\Del_{\m} $ is the highest weight vector of
the weighted Bergman space $L^2_a(\Ome, S_{\m}(V), \mu_{\a})$,
 and $\bar D^m$
intertwines the $G$-action $\pi_{\nu}$ on $A^{2, \a}_{\m}(\Ome)$
with that on  $L^2_a(\Ome, S_{\m}(V), \mu_{\a})$ (see (4.2)), by Lemma 2.1. Thus  it is
a non-zero intertwining operator of the two spaces.
We summarize our results  in the following
\begin{theo+}\label{main-th} The relative discrete series $A^{2, \alpha}_{\m}(\Ome)$ 
is $G$-equivalent to the weighted Bergman space $L^2_a(\Ome, V^{\m}, \mu_{\a})
$ and the corresponding intertwining operator is given by
$\bar D^{|\m|}$. The highest weight vector of  $A^{2, \alpha}_{\m}(\Ome)$ 
is given by $\Del_{\m}(q(z))$. In particular, the space $A^{2, \alpha}_{\m}(\Ome)$ 
consists of nearly holomorphic functions.
\end{theo+}

\begin{rema+} By the results of Shimeno \cite{Shimeno} we see that all
the relative discrete series are obtained in this way.
\end{rema+}
\begin{rema+} When $\Ome$ is of tube type and when $\m
=m\ub1r=(m, m, \dots, m)$, the above
is also proved in \cite{ahd+bo+gkz} by considering
the tensor products of Bergman spaces
holomorphic functions with polynomial space
of anti-holomorphic functions,   and 
in \cite{oz-reldis} by Capelli identity.
\end{rema+}

\begin{section}{An example: The case of the unit ball in
$\mathbb C^n$}
\end{section}

In this section we consider the example of the unit ball
in $V=\mathbb C^n$. 
The rank  $r=1$, $h(z)=1-           |z|^2$, 
 and symmetric tensor 
$S_{m}(V)$  is  itself
irreducible under $K$.
We study the
adjoint operator $\bar D^\ast$ 
of $\bar D$
instead of $\bar D$. 
The operator $(\bar D^\ast)^m$
is thus an intertwining operator from the weighted Bergman space
$L^2_a(\Ome, S_{m}(V), \mu_{\a})$
of vector-valued holomorphic functions
 into
the relative discrete series $A^{2, \a}_{m}(\Ome)$.
Now the $L^2_a(\Ome, S_{m}(V), \mu_{\a})$
has highest weight vector $\otimes^m e_1$. Thus
$(\bar D^\ast)^m (\otimes^m e_1)$
is the highest weight vector in $A^{2, \a}_{m}$. We calculate
directly here this vector.

Let $D=\bar D^\ast $. It has the following expression on
a function $f$ with values in $\otimes^m V$:
$$
Df=h(z)^{-\a} \otimes^{m-1}B(z, \bar z) \Tr \partial\,
\left[h(z)^{\a}(I\otimes \otimes^{m-1}B(z,\bar  z)^{-1})f\right].
$$
To explain the formula we note that, the operator $\partial$
acting on a $\otimes^m V$-valued function gives a
functions with values in $V^\prime\otimes(\otimes^m V)
=(V^\prime\otimes V)
\otimes(\otimes^{m-1} V)$;
the operator $\Tr$ is the bilinear pairing between
the first  factor $V^\prime\otimes V$.
Recall that  the Bergman operator on the unit ball is
$$
B(z, \bar z)=(1-|z|^2)(1-z\otimes z^\ast),
$$
where $z\otimes z^\ast$ is the rank one operator on 
$V$, $z\otimes z^\ast(v)=\langle v, z\rangle z$; see \cite{pz-cr}.
Take $f =\otimes^m e_1$. The above formula then reads
\begin{equation*}
\begin{split}
&\quad \, D\otimes^m e_1\\
& =
h(z)^{-\a-(m-1)} \otimes^{m-1}(1-z\otimes z^\ast) 
\Tr \partial
\left[h(z)^{\a-2(m-1)}(e_1\otimes \otimes^{m-1}((1-|z|^2)e_1+\bar z_1z))\right]
\end{split}
\end{equation*}
Performing the differentiation using the Leibniz rule we first
differentiate the term $h(z)^{\a-2(m-1)}$, and get
\begin{equation}
(2(m-1)-\a)
h(z)^{\a-2(m-1)-1}(\sum_j \bar z_j dz_j)\otimes
(e_1\otimes \otimes^{m-1}((1-|z|^2)e_1+\bar z_1z).
\end{equation}
Taking the trace $\Tr$, it is 
$$
(2(m-1)-\a) (1-|z|^2)^{-1}\bar z_1\otimes^{m-1} e_1.
$$ 
Next we differentiate each factor
 $(1-|z|^2)e_1+\bar z_1z$ in the tensor, and get
\begin{equation}
(-\sum_j \bar z_j dz_j)\otimes e_1+\bar z_1
\sum_j dz_j\otimes e_j. 
\end{equation}
We perform the operation $Tr$  and observe that
each term  is vanishing:
\begin{equation}\Tr e_1\otimes
((-\sum_j \bar z_j dz_j)\otimes e_1+\bar z_1
\sum_j dz_j\otimes e_j)
=0.
\end{equation}
Thus only the first differentiation contributes 
to the final result, that is
$$
D (\otimes^m e_1)
=(2(m-1)-\a) (1-|z|^2)^{-1}\bar z_1\otimes^{m-1} e_1.
$$
By induction we get 
$$
D^m e_1^m
=C(1-|z|^2)^{-m}\bar z_1^m
$$
where{\sl \/}
$$
C=\prod_{l=0}^{{m-1}}(2(m-1-j)-\a +j).
$$

The function $(1-|z|^2)^{-m}\bar z_1^m$ is in 
$L^2(\Ome, \mu_{\a})$ if and only if $0\le m<\frac {\a+1}2$. In that case
 $D^m (\otimes^m e_1)$ is a non-zero multiple of $(1-|z|^2)^{-m}\bar z_1^m$.
The quasi-inverse  is 
$q(z)={(1-|z|^2)^{-1}}{\bar z}
$, and
the vector constructed in Theorem \ref{main-th}
is $[e_1^m, \otimes q(z)]=(1-|z|^2)^{-m}\bar z_1^m$,
and thus the two methods give the same result.

One might also in the beginning work with the operator
$D=-(\bar D)^\ast$  instead of $\bar D$.
However we note that
for a general bounded symmetric domain
the formula for the operator $D$
is much more involved.

\end{document}